\documentclass[12pt]{article}
\usepackage{graphicx, mathrsfs,amssymb,amsmath,amsfonts,amsbsy,amsthm,latexsym,bm,float}
\begin{document}

\newcommand{\Dh}{\hbox{\bf D}}
\newcommand{\Xh}{\hbox{\bf X}}
\newcommand{\Fh}{\hbox{\bf F}}
\newcommand{\nablaslh}{\nabla\raise2pt\hbox{\kern-9pt\slash\kern3pt}}

\title{
Quadrature formulas for integrals transforms generated by orthogonal polynomials\\
\vskip.5cm}
\author{Rafael G. Campos, Francisco Dom\'{\i}nguez Mota,  E. Coronado\\
Facultad de Ciencias F\'{\i}sico-Matem\'aticas,\\
Universidad Michoacana, \\
58060, Morelia, Mich., M\'exico.\\
\hbox{\small rcampos@umich.mx, dmota@umich.mx, ecoronado@fismat.umich.mx}\\
}
\date{}
\maketitle
{
\noindent MSC: 33C45, 33C47, 44A20, 65D32\\
\noindent Keywords: integral transforms, quadrature, orthogonal polynomials, bilinear generating functions.
}\\
\vspace*{1truecm}
\begin{center} Abstract \end{center}
By using the three-term recurrence equation satisfied by a family of orthogonal polynomials, the Christoffel-Darboux-type bilinear generating function and their asymptotic expressions, we obtain quadrature formulas for integral transforms generated by the classical orthogonal polynomials. These integral transforms, related to the so-called Poisson integrals, correspond to a modified Fourier Transform in the case of the Hermite polynomials, a Bessel Transform in the case of the Laguerre polynomials and to an Appell Transform in the case of the Jacobi polynomials.
\vskip1.5cm
\newpage
\section{Introduction}\label{intro}
Certain integral transforms with a Christoffel-Darboux-type bilinear generating function for a classical orthogonal polynomial in the kernel, termed Poisson integrals \cite{Muc69, Erd41}, appear in the expansion of functions in terms of the classical orthogonal polynomials and, under a suitable change of variable, they become well-known integral transforms for certain limit values of the expansion parameter. In this sense, Mehler's formula yields a modified Fourier transform in the case of the Hermite polynomials, the Hille-Hardy formula gives a Bessel-Hankel transform in the case of the Laguerre polynomials and Bailey's bilinear generating function produces an Appell transform in the case of the Jacobi polynomials.\\
We show in this paper that each of these integral transforms has a quadrature formula of gaussian type. This was done by following {\sl mutatis mutandis} the approach of \cite{Cam92, Cam95}, where a quadrature formula for the Fourier transform and other one for the Hankel transform were obtained by using the recurrence equation for the (generalized) Hermite polynomials and their zeros as well as the Christoffel-Darboux formula, Mehler's formula and the Hille-Hardy formula. The generalization of this procedure provides a method to yield at one time new integral transforms and their quadrature formulas. In the case of a family of orthogonal polynomials ${\mathscr P}=\{P_k(x), k=0,1,\ldots\}$ this method requires a suitable bilinear generating function of ${\mathscr P}$ as well as some knowledge about the asymptotic behavior of the polynomials and their zeros. The kernel of the integral transform becomes related to the generating function. We apply this technique to the classical orthogonal polynomials to obtain quadrature formulas for their corresponding Poisson integrals. Thus, in Sec. \ref{sectre} we present a quadrature for a modified Fourier transform generated by the Hermite polynomials, a quadrature for a modified Bessel transform yielded by the Laguerre polynomials in Sec. \ref{seccua}, and finally, in Sec. \ref{seccin} we introduce a quadrature for an Appell transform generated by the Jacobi polynomials. This integral transform has also been obtained in \cite{Vir92}.
\section{Outline of the method}\label{secdos}
Let ${\mathscr P}=\{P_k(x), k=0,1,\ldots\}$ be one of the families of classical orthogonal polynomials (Hermite, Laguerre or Jacobi) satisfying the recurrence equation
\begin{equation}\label{receqg}
A_nP_{n+1}(x)+B_nP_n(x)+C_{n-1}P_{n-1}(x)=xP_n(x), \quad n=0,1,2,\ldots,
\end{equation}
where $ P_{-1}(x)\equiv 0$. Let $k_n$ denote the coefficient of $x^n$ in $P_n(x)$. As it is well-known \cite{Sze75, Chi78}, from (\ref{receqg}) follows the Christoffel-Darboux formula 
\begin{equation}\label{chrsdar}
\sum_{n=0}^{N-1} s_n^2 P_n(x)P_n(y)=\frac{k_{N-1}s^2_{N-1}}{k_N}\,\displaystyle\frac{P_N(x)P_{N-1}(y)-P_{N-1}(x)P_N(y)}{x-y},
\end{equation}
where $s_n$ is determined by the reciprocal of the norm 
\begin{equation}\label{dessk}
\left(\int_a^bP_n^2(x)d\omega(x)\right)^{1/2}
\end{equation}
up to a numerical constant independent of $n$. 
Here, $(a,b)$ is the orthogonality interval of ${\mathscr P}$ and $d\omega(x)$ is the corresponding non-negative measure. The recurrence equation (\ref{receqg}) can be written as the eigenvalue problem $M_\infty P= x P$, where $(M_\infty)_{nk}=A_{n-1}\delta_{n+1,k}+B_{n-1}\delta_{n,k}+C_{n-2}\delta_{n,k+1}$, $n,k=1,2,\ldots$, and $P$ is the vector whose $n$th entry is $P_{n-1}(x)$. Let us now consider the eigenproblem associated to the principal submatrix of dimension $N$ of $M_\infty$. This part is a well-known technique \cite{Wil62, Gol69, Gau99} to yield gaussian quadratures: the three-term recurrence equation is rewritten in matrix form to obtain orthonormal vectors of ${\mathbb R}^N$ whose entries are given in terms of the values of $P_k(x)$, $k=0,1,\ldots,N-1$, at the zeros of $P_N(x)$. To proceed, we take a similarity transformation to symmetrize $M_\infty$. The diagonal matrix $S=\text{diag}\{s_0,s_1,\ldots,\}$ whose elements $s_n$ are given by (\ref{dessk}), generates a symmetric matrix $SM_\infty S^{-1}$ whose principal submatrix of order $N$, denoted by $M$, has elements given by 
\[
M_{nk}=\sqrt{A_{n-1}C_{n-1}} \delta_{n+1,k}+B_{n-1}\delta_{n,k}+\sqrt{A_{n-2}C_{n-2}} \delta_{n,k+1}, \quad n,k=1,2,\ldots,N.
\]
The recurrence equation (\ref{receqg}) and formula (\ref{chrsdar}) can be used to solve the eigenproblem 
\[
Mu_k=x_k u_k,\quad k=1,2,\ldots, N.
\]
The eigenvalues $x_k$ are the zeros of $P_N(x)$ and the $n$th entry of the $k$th eigenvector $u_k$ is given by $c_k s_{n-1}P_{n-1}(x_k)$, $n=1,2,\ldots,N$, where $c_k$ is a normalization constant which is obtained from (\ref{chrsdar}). Thus, the orthonormal vectors $u_k$, $k=1,2,\ldots, N$, have components
\begin {equation}\label{ortvec}
u_{n-1}(x_k)=\left(\displaystyle\frac{s_{n-1}}{s_{N-1}}\right)
\frac{P_{n-1}(x_k)}{\sqrt{k_{N-1}P_{N-1}(x_k)P'_N(x_k)/k_N}},\quad n=1,\ldots,N.
\end{equation}
Note that the product $k_{N-1}P_{N-1}(x_k)P'_N(x_k)/k_N$ is always positive since it is $1/s^2_{N-1}$ times the squared norm of the vector 
\[
\left(s_0P_0(x_k),s_1P_1(x_k),\ldots,s_{N-1}P_{N-1}(x_k)\right).
\]
Let $U$ be the orthogonal matrix whose $k$th column is $u_k$ and $D$ be the diagonal matrix $D=\text{diag}\{1,z,z^2,\ldots,z^{N-1}\}$, where $z$ is a complex number. Then, the matrix $T=U^{-1}DU$ whose elements are explicitly given by
\begin{eqnarray}\label{tmat}
T_{jk}&=&\sum_{n=0}^{N-1}z^nu_n(x_j)u_n(x_k)\nonumber\\
&=&\displaystyle\frac{\vert k_N/k_{N-1}\vert}{s^2_{N-1}\sqrt{\vert P_{N-1}(x_j)P'_N(x_j)P_{N-1}(x_k)P'_N(x_k)\vert}}
\sum_{n=0}^{N-1}s^2_nP_n(x_j)P_n(x_k)z^n,
\end{eqnarray}
is the Discrete Transform associated to the corresponding Integral Transform. To show this, we use some asymptotic properties of the zeros of the classical orthogonal polynomials shown for each case in the next sections. 
First, we note that the asymptotic expressions for $P_{N-1}(x)$ and $P'_N(x)$ in the oscillatory region, evaluated at the zeros of $P_N(x)$ satisfy the formula 
\begin{equation}\label{prodpols}
\sqrt{\vert P_{N-1}(x_k)P'_N(x_k)\vert}= \frac{a_N}{g(x_k)}+ O(N^{-\mu}).
\end{equation}
The function $g(x)$ is positive in $(a,b)$ and it is related to the weight function of ${\mathscr P}$. The constant $a_N>0$ is independent of $k$ and the exponent $\mu$ can be taken positive. The zeros of $P_N(x)$ contained in the fixed interval $(c,d)$, $a<c<d<b$, become evenly spaced for large values of $N$ under a bijective mapping $y=\sigma(x)$, i.e.,
\begin{equation}\label{Dely}
\Delta y_k=\Delta \sigma(x_k)=\sigma(x_{k+1})-\sigma(x_k)=\lambda_N+ O(N^{-\mu}),
\end{equation}
where $\lambda_N$ does not depend on $k$. Furthermore, we have that
\begin{equation}\label{rmed}
\displaystyle\frac{\vert k_N/k_{N-1}\vert}{s^2_{N-1}a^2_N}= A\Delta \sigma(x_k)+ O(N^{-\mu}),
\end{equation}
where $A$ is a numerical constant independent of $N$ and $k$. Since the Christoffel-Darboux-type bilinear generating function $G(x,y,z)$ of the set ${\mathscr P}$ satisfies
\begin{equation}\label{bilgenfun}
\sum_{n=0}^\infty s^2_n P_n(x)P_n(y)z^n=G(x,y,z),
\end{equation}
for $(x,y)\in (a,b)\times (a,b)$ and $z$ in a suitable domain ${\mathcal D}$ of the complex plane, the element $T_{jk}$ [cf. Eq. (\ref{tmat})] has the limiting form
\begin{equation}\label{cuadprev}
A\,g(x_j)g(x_k)G(x_j,x_k,z)\Delta \sigma(x_k),
\end{equation}
for large values of $N$. Therefore, 
\begin{equation}\label{tmatasymdos}
\sum_{k=1}^NT_{jk}f(x_k)\to A \sum_{k=1}^Ng(x_j)g(x_k)G(x_j,x_k,z)f(x_k)\Delta \sigma(x_k), \quad N\to\infty.
\end{equation}
The sum of the right-hand side of this equation is the Riemann-Stieltjes sum of $f(x)$ with respect to $\sigma(x)$ and tends to the integral transform
\begin{equation}\label{intrasnf}
{\mathcal T}[f(x);y,z]=\int_a^b  K(x,y,z) f(x)d\sigma(x),
\end{equation}
where $K(x,y,z)=Ag(x)g(y)G(x,y,z)$. Therefore, Eq. (\ref{tmatasymdos}) becomes the quadrature formula
\begin{equation}\label{quadforgen}
\int_a^b K(x,y_j,z) f(x)d\sigma(x)\simeq\sum_{k=1}^NT_{jk}f(x_k),
\end{equation}
for $z\in{\mathcal D}$.
\section{A modified Fourier transform}\label{sectre}
As it is known, in the case in which $P_n(x)$ is the $n$th Hermite polynomial $H_n(x)$, the parameters $k_n$ and $s^2_n$ are $2^n$ and $1/(2^n n!)$ respectively and the nodes $x_k$ are zeros of $H_N(x)$. Taking into account that $\vert H'_N(x_k)\vert=(-1)^{N+k}H_{N-1}(x_k)$, the $n$th component of the $k$th orthonormal vector $u_k$ given by Eq. (\ref{ortvec}) becomes
\begin {equation}\label{orther}
u_{n-1}(x_k)=(-1)^{N+k}\left(\displaystyle\frac{2^{N-n}(N-1)!}{N(n-1)!}\right)^{1/2}\frac{H_{n-1}(x_k)}{H_{N-1}(x_k)},\quad n=1,\ldots,N.
\end{equation}
Therefore, the elements of the Discrete Transform $T$, denoted in this case by $T_H$, are
\begin {equation}\label{Ther}
(T_H)_{jk}=(-1)^{j+k}\frac{2^{N-1}(N-1)!}{NH_{N-1}(x_j)H_{N-1}(x_k)}\sum_{n=0}^{N-1}\frac{z^n}{2^nn!}H_n(x_j)H_n(x_k).
\end{equation}
By using the asymptotic formula \cite{Sze75}
\[
H_N(x)=\frac{\Gamma(N+1)}{\Gamma(N/2+1)}e^{x^2/2}\left(\cos(\sqrt{2N+1}\,\,x-N\pi/2)+{\cal O}(N^{-1/2})\right),
\]
we obtain that $g(x)=e^{-x^2/2}$, $\sigma(x)=x$, $\Delta x_k=\pi/\sqrt{2N}$, $A=1$ [cf. Eqs (\ref{prodpols})-(\ref{rmed})] and the expression for $G(x,y,z)$ is determined by Mehler's formula. Thus, we get the quadrature formula
\begin{equation}\label{herquad}
\int_{-\infty}^\infty K_H(x,y_j,z) f(x)dx\simeq\sum_{k=1}^N(T_H)_{jk}f(x_k),
\end{equation}
where
\[
K_H(x,y,z)=\frac{1}{\sqrt{\pi(1-z^2)}}\,\exp\left(-\frac{(1+z^2)(x^2+y^2)-4xyz}{2(1-z^2)}\right),\quad \vert z\vert<1.
\]
The argument of the exponential is imaginary if $z\ne\pm 1$ is on the unit circle. In particular, if $z=\pm i$, $K_H(x,y,z)=\exp(\pm ixy)/\sqrt{2\pi}$, and (\ref{herquad}) becomes a quadrature for the Fourier Transform. This formula has been previously obtained in \cite{Cam92} were some numerical examples also has been given. 
\subsection{An example}
\begin{figure}[H]\label{Figuna}
\centering
\includegraphics[scale=1.0]{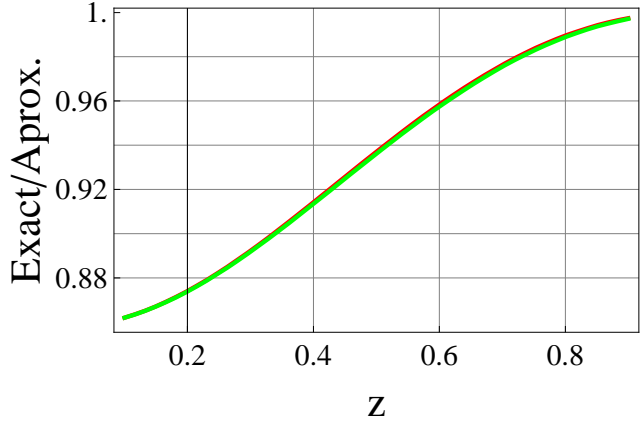}
\caption{Real part of the left-hand side against the real part of the right-hand side of (\ref{fqexuno}) for $N=31$ and $z\in (0,1)$. The norm of the error is 0.0048.}
\end{figure}

We fix the point $y_j$ where (\ref{herquad}) is evaluated and thus, the integral becomes a function only of $z$. For $z$ real, the kernel is a real exponential but even so the Fourier transform can be recovered by taking $f(x)$ as a complex exponential. In this example we take $y_j=0$ ($N$ odd) and $f(x)=\exp(-ix)$. Thus, (\ref{herquad}) becomes
\begin{equation}\label{fqexuno}
\frac{2}{1+z^2} e^{-\frac{1-z^2}{2(1+z^2)}}\simeq \sum_{k=1}^N(T_H)_{\frac{N+1}{2},k}\,e^{-ix_k}.
\end{equation}
In Fig. \ref{Figuna} we show the numerical computation of this formula.
\section{A Bessel transform }\label{seccua}
We take now the Laguerre polynomial $L_n^{(\alpha)}(x)$ as $P_n(x)$. Thus, $k_n=(-1)^n/n!$, $s^2_n=n!/\Gamma(\alpha+n+1)$ and the nodes $x_k$ are the $N$ zeros of $L_N^{(\alpha)}(x)$. We have now that $\vert {L'}_N^{(\alpha)}(x_k)\vert=(-1)^{k+1}\sqrt{(N+\alpha)/x_k}L_{N-1}^{(\alpha)}(x_k)$ and the vector $u_k$ given by (\ref{ortvec}) becomes
\begin {equation}\label{ortlag}
u_{n-1}(x_k)=(-1)^{k+1}\left(\displaystyle\frac{(n-1)!\Gamma(\alpha+N) x_k}{N!(N+\alpha)\Gamma(\alpha+n)}\right)^{1/2}\frac{L_{n-1}^{(\alpha)}(x_k)}{L_{N-1}^{(\alpha)}(x_k)},\quad n=1,\ldots,N.
\end{equation}
The matrix $T$, denoted now by $T_L$ has components
\begin {equation}\label{Tlag}
(T_L)_{jk}=\frac{(-1)^{j+k}\Gamma(\alpha+N) \sqrt{x_jx_k}}{N!(N+\alpha)L_{N-1}^{(\alpha)}(x_j)L_{N-1}^{(\alpha)}(x_k)}\sum_{n=0}^{N-1}\frac{n!\,z^n}{\Gamma(\alpha+n+1)}L_n^{(\alpha)}(x_j)L_n^{(\alpha)}(x_k).
\end{equation}
The use of the asymptotic formula \cite{Sze75}
\[
L_N^{(\alpha)}(x)=\frac{N^{\alpha/2-1/4}e^{x/2}}{\sqrt{\pi}x^{\alpha/2+1/4}}\left(\cos(2\sqrt{Nx}-\alpha\pi/2-\pi/4)+(Nx)^{-1/2}{\cal O}(1)\right)
\]
yields $g(x)=x^{\alpha/2+1/4}e^{-x/2}$, $\sigma(x)=\sqrt{x}$, $\Delta\sigma(x_k)=\pi/2\sqrt{N}$, $A=1/2$. The generating function $G(x,y,z)$ is given by the Hille-Hardy formula. Thus, we get the quadrature formula
\begin{equation}\label{lagquad}
\int_0^\infty K_L(x,y_j,z) f(x)dx\simeq\sum_{k=1}^N(T_L)_{jk}f(x_k),
\end{equation}
where
\[
K_L(x,y,z)=\frac{z^{-\alpha/2}}{1-z}\left(\frac{y}{x}\right)^{1/4}\exp\left(-\frac{(1+z)(x+y)}{2(1-z)}\right)I_\alpha\left(\frac{2\sqrt{xyz}}{1-z}
\right),\quad \vert z\vert< 1,
\]
and $I_\alpha(x)$ is the modified Bessel function of the first kind. Note that we have used the fact that $d\sigma(x)=x^{-1/2}dx/2$.\\
For $z=-1$, (\ref{lagquad}) becomes a quadrature for the Hankel transform
\begin{equation}\label{hanquad}
\frac{1}{2}\,y_j^{1/4}\int_0^\infty x^{-1/4} J_\alpha(\sqrt{y_jx})  f(x)dx\simeq \sum_{k=1}^N(T_L)_{jk}f(x_k),
\end{equation}
which is written in an nonstandard way.
\subsection{An example}
\begin{figure}[H]\label{Figdos}
\centering
\includegraphics[scale=1.0]{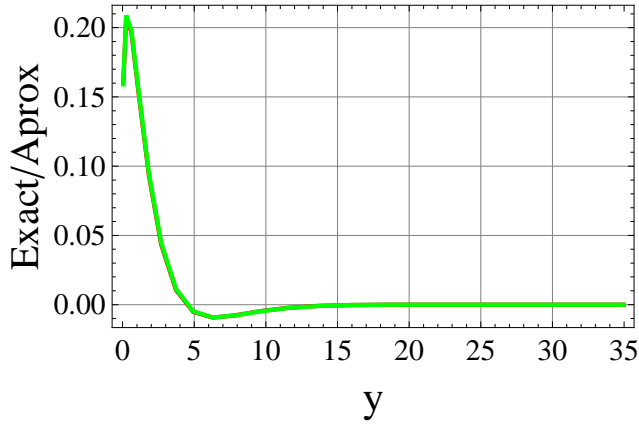}
\caption{The left-hand side of (\ref{fqexdos}) is plotted versus the right-hand side.  The values of the parameters are $N=30$, $\alpha=0$, $c=2$, and $z=1/10$. The norm of the error is 0.0036.}
\end{figure}

The integral
\begin{equation}\label{exados}
\int_{0}^\infty e^{-a x^2} x J_\alpha(cx)I_\alpha(bx) dx=\frac{1}{2a}e^{\frac{b^2-c^2}{4a}}J_\alpha(\frac{cb}{2a})
\end{equation}
is given usually as a Hankel transform \cite{Erd53b}, but it will be used here to test (\ref{lagquad}). To this end, we make the change of variable $x\to x^2$ in the left-hand side of (\ref{lagquad}) and take 
\[
a=\frac{1+z}{2(1-z)},\quad b_j=\frac{2\sqrt{zy_j}}{1-z},\quad f(x^2)=\sqrt{x}J_\alpha(cx).
\]
Thus, (\ref{lagquad}) becomes
\begin{equation}\label{fqexdos}
\frac{2 z^{-\alpha/2}}{1+z}y_j^{1/4}\exp\left[-\left(\frac{1-z}{1+z}\right)\frac{c^2+y_j}{2}\right] J_\alpha\left(\frac{2c\sqrt{y_j z}}{1+z}\right)\simeq
\sum_{k=1}^N(T_L)_{jk}x_k^{1/4}J_\alpha(c\sqrt{x_k}).
\end{equation}
The numerical output of this quadrature is illustrated in Fig. \ref{Figdos}. Note that (\ref{fqexdos}) is an {\it exact} formula if $z=1$ since $T$ becomes the identity matrix for this value of $z$.
\section{An Appel transform }\label{seccin}
Proceeding as before, we take now the Jacobi polynomial $P_n^{(\alpha,\beta)}(x)$ as $P_n(x)$. Therefore, 
\[
k_n=\frac{(2n+\alpha+\beta)!}{2^n (n+\alpha+\beta)n!},\qquad 
s^2_n=\frac{(2n+\alpha+\beta+1)n!\Gamma(n+\alpha+\beta+1)}{\Gamma(n+\alpha+1)\Gamma(n+\beta+1)}.
\]
In this case, the nodes $x_k$ are the zeros of $P_N^{(\alpha,\beta)}(x)$. We have now that 
\[
\vert {P'}_N^{(\alpha,\beta)}(x_k)\vert=(-1)^{N+k}\frac{2(N+\alpha)(N+\beta)}{(2N+\alpha+\beta)(1-x_k^2)}P_{N-1}^{(\alpha,\beta)}(x_k)
\]
and the components of the $k$th orthonormal vector $u_k$ are
\begin {eqnarray}
u_{n-1}(x_k)&=&(-1)^{N+k}\bigg(\displaystyle\frac{(n-1)!(2n+\alpha+\beta-1)(2N+\alpha+\beta)^2\Gamma(N+\alpha)\Gamma(N+\beta)}
{4N!(N+\alpha)(N+\beta)\Gamma(N+\alpha+\beta+1)\Gamma(n+\alpha)}\nonumber\\
&\times&\frac{\Gamma(n+\alpha+\beta)}{\Gamma(n+\beta)}\bigg)^{1/2}
\sqrt{1-x_k^2}\,\,\frac{P_{n-1}^{(\alpha,\beta)}(x_k)}{P_{N-1}^{(\alpha,\beta)}(x_k)},\quad n=1,\ldots,N.\label{ortjac}
\end{eqnarray}
Therefore, the components of the Discrete Transform $T$, denoted here by $T_J$ are
\begin {eqnarray}
(T_J)_{jk}&=&(-1)^{j+k}\displaystyle\frac{(2N+\alpha+\beta)^2\Gamma(N+\alpha)\Gamma(N+\beta)}
{4N!(N+\alpha)(N+\beta)\Gamma(N+\alpha+\beta+1)}\frac{\sqrt{(1-x_j^2)(1-x_k^2)}}{P_{N-1}^{(\alpha,\beta)}(x_j)P_{N-1}^{(\alpha,\beta)}(x_k)}
\nonumber\\
&\times&\sum_{n=0}^{N-1}\frac{n!(2n+\alpha+\beta+1)\Gamma(n+\alpha+\beta+1)}{\Gamma(n+\alpha+1)\Gamma(n+\beta+1)}\,z^n
P_n^{(\alpha,\beta)}(x_j)P_n^{(\alpha,\beta)}(x_k).\label{Tjac}
\end{eqnarray}
The use of the asymptotic formula \cite{Sze75}
\[
P_N^{(\alpha,\beta)}(\cos\theta)=\frac{1}{\sqrt{N\pi}}\frac{\cos\left[\left(N+\frac{\alpha+\beta+1}{2}\right)\theta-\left(\alpha+\frac{1}{2}\right)\frac{\pi}{2}\right]}
{\left(\sin\frac{\theta}{2}\right)^{\alpha+1/2}\left(\cos\frac{\theta}{2}\right)^{\beta+1/2}}+{\cal O}(N^{-3/2})
\]
yields $g(x)=(1-x)^{\alpha/2+1/4}(1+x)^{\beta/2+1/4}$, $\sigma(x)=\arccos(x)$, $\Delta\sigma(x_k)=\pi/N$, $A=2^{-\alpha-\beta-1}$. The generating function $G(x,y,z)$ is given by Bailey's formula \cite{Bai38}
\begin{eqnarray*}
G_J(x,y,z)&=&\frac{\Gamma(\alpha+\beta+2)}{\Gamma(\alpha+1)\Gamma(\beta+1)}\frac{(1-z)}{(1+z)^{\alpha+\beta+2}}\\
&\times& F_4\left[\frac{\alpha+\beta+2}{2},\frac{\alpha+\beta+3}{2};\alpha+1,\beta+1;\frac{z(1-x)(1-y)}{(1+z)^2},\frac{z(1+x)(1+y)}{(1+z)^2}\right],
\end{eqnarray*}
where $F_4(a,b;c,d;\xi,\eta)$ is the fourth Appel's hypergeometric function of two variables \cite{Erd53}. If $r$ and $s$ are the radii of convergence in $\xi$ and $\eta$ respectively, then $r^{1/2}+s^{1/2}=1$ and $x$, $y$ and $z$ should satisfy
\[
\left\vert\frac{z(1-x)(1-y)}{(1+z)^2}\right\vert < r, \quad \left\vert\frac{z(1+x)(1+y)}{(1+z)^2}\right\vert< s,
\]
in order to have the quadrature formula
\begin{equation}\label{jacquad}
\int_{-1}^1 K_J(x,y_j,z) f(x)dx\simeq\sum_{k=1}^N(T_J)_{jk}f(x_k),
\end{equation}
where
\[
K_J(x,y,z)=\frac{[(1-x)(1-y)]^{\alpha/2+1/4}[(1+x)(1+y)]^{\beta/2+1/4}}{2^{\alpha+\beta+1}\sqrt{1-x^2}}\,\,G_J(x,y,z).
\]
Note that we have written the integral in terms of $dx$ instead of $d\sigma(x)$.
\subsection{An example}
\begin{figure}[H]\label{Figtres}
\centering
\includegraphics[scale=1.0]{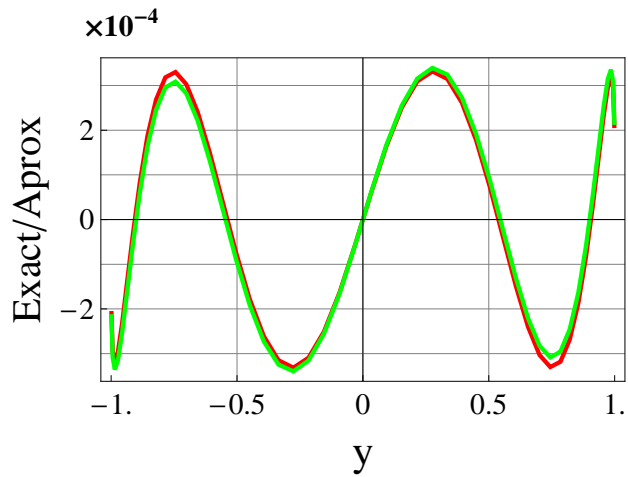}
\caption{Plot of the left-hand side against the right-hand side of (\ref{fqextres}) for $N=50$, $\alpha=\beta=0$, $z=1/4$ and $n=5$. The norm of the error is  0.00011 in this case.}
\end{figure}

It is shown in \cite{Vir92} that\footnote{We have corrected a misprint in this formula.}
\begin{eqnarray}\label{exatres}
\int_{-1}^1& (1-x)^\alpha(1+x)^\beta P_n^{(\alpha,\beta)}(x) F_4(\bar{\alpha},\bar{\beta};\alpha+1,\beta+1;\xi,\eta) dx\nonumber\\
&=\displaystyle\frac{2^{\alpha+\beta+1}\Gamma(\alpha+1)\Gamma(\beta+1)}{(\alpha+\beta+1)\Gamma(\alpha+\beta+1)}\frac{(1+z)^{\alpha+\beta+2}}{(1-z)}z^nP_n^{(\alpha,\beta)}(y),
\end{eqnarray}
where
\[
\bar{\alpha}=\frac{\alpha+\beta+2}{2},\quad \bar{\beta}=\frac{\alpha+\beta+3}{2},\quad \xi=\frac{z(1-x)(1-y)}{(1+z)^2},\quad \eta=\frac{z(1+x)(1+y)}{(1+z)^2}.
\]
By choosing a suitable integrand, our quadrature formula (\ref{jacquad}) becomes 
\begin{equation}\label{fqextres}
z^n (1-y_j)^{\alpha/2+1/4}(1+y_j)^{\beta/2+1/4} P_n^{(\alpha,\beta)}(y_j)\simeq\sum_{k=1}^N(T_J)_{jk} (1-x_k)^{\alpha/2+1/4}(1+x_k)^{\beta/2+1/4} P_n^{(\alpha,\beta)}(x_k)
\end{equation}
for this case. The numerical output of this quadrature is illustrated in Fig. \ref{Figtres}. This equation is an {\it exact} formula if $z=1$.
\section{Final Remarks}\label{secseis}
There are important facts supporting this method for finding quadrature formulas for the integral transforms (Poisson integrals) associated to a family of orthogonal polynomials. These facts are given by formulas (\ref{prodpols})-(\ref{bilgenfun}). An important one, is that the kernel of the integral transform is determined by the product of the bilinear generating function $G(x,y,z)$ and $g(x)g(y)$. The bilinearity of the kernel explains why the associated quadrature formula can yield exact results if the function to be transformed is chosen properly in terms of elements of the orthogonal family.

\end{document}